\documentclass{opt2025_v2} 

\usepackage{amsmath,amssymb,amsfonts}%
\usepackage{graphicx}%
\usepackage{multirow}%
\usepackage{mathrsfs}%
\usepackage{xcolor}%
\usepackage{textcomp}%
\usepackage{manyfoot}%
\usepackage{booktabs}%
\usepackage{algorithm}%
\usepackage{algorithmicx}%
\usepackage{algpseudocode}%
\usepackage{listings}%
\newtheorem{assumption}[theorem]{Assumption}
\usepackage{subcaption}
\usepackage{comment}
\usepackage{todonotes}
\newcommand{\ze}[1]{{\color{black}#1}}
\newcommand{\aj}[1]{{\color{black}#1}}
\newcommand{\z}[1]{{\color{black}#1}}
\newcommand{\za}[1]{{\color{black}#1}}
\newcommand{\zz}[1]{{\color{black}#1}}
\newcommand{\af}[1]{{\color{black}#1}}
\title[Distributionally Robust Nash Games]{Distributionally Robust Nash Equilibria via Variational Inequalities}

\optauthor{%
\Name{Zeinab Alizadeh} \Email{zalizadeh@arizona.edu}\\
\Name{Azadeh Farsi} \Email{azadehfarsi@arizona.edu}\\
\Name{Afrooz Jalilzadeh} \Email{afrooz@arizona.edu}\\
\addr Department of Systems and Industrial Engineering, The University of Arizona, USA}

\allowdisplaybreaks
\begin{document}
\allowdisplaybreaks
\maketitle

\begin{abstract}%
Nash Equilibrium and its robust counterpart, Distributionally Robust Nash Equilibrium (DRNE), are fundamental problems in game theory with applications in economics, engineering, and machine learning. This paper addresses the problem of DRNE, where multiple players engage in a noncooperative game under uncertainty. Each player aims to minimize their objective against the worst-case distribution within an ambiguity set, resulting in a minimax structure. We reformulate the DRNE problem as a Variational Inequality (VI) problem, providing a unified framework for analysis and algorithm development. We propose a gradient descent-ascent type algorithm with convergence guarantee that effectively addresses the computational challenges of high-dimensional and nonsmooth objectives. 
\end{abstract}

\section{INTRODUCTION}\label{sec:intro}

\aj{Noncooperative game theory provides a rigorous mathematical framework for analyzing multi-agent decision-making problems that arise in a wide range of applications in economics, engineering, and machine learning. Within this domain, Nash games are a fundamental modeling tool, where a finite set of self-interested agents competes, each seeking to optimize their own objectives while accounting for the strategies of others.
However, such systems often involve uncertainty, where key parameters—such as payoffs or costs—are not precisely known. This uncertainty necessitates robust formulations that can handle ambiguity and ensure reliable solutions, even under worst-case conditions.
To be more specific, consider the Stochastic Nash Equilibrium (SNE) problem:
\begin{align}\label{SNE}\tag{SNE}
\min _{x_i \in K_i} \mathbb{E}\left[f_i\left(x_i, x_{-i}, \xi_i\right)\right], \quad i = 1, \ldots, n,
\end{align}
where $x_i$ denotes the strategy of the $i$th player, $x_{-i}$ is the collection of the strategies of the other players, $\xi_i$ is a random vector whose probability distribution is unknown but belongs to an ambiguity set, $f_i$ is a convex and possibly nonsmooth function, and $K_i\subseteq \mathbb R^{n_i}$ is a convex and compact set. The expectation operator $\mathbb{E}[\cdot]$ represents the mathematical expectation over the uncertain parameters. 

When the probability distribution of $\xi_i$ is not fully specified, a distributionally robust formulation provides a more reliable alternative. Specifically, we reformulate the game as:
\begin{equation}
\label{DRO_Nash}\tag{DRNE}
\min _{x_i \in K_i} \max _{p_i \in P_i} \sum_{j=1}^m p_{ij} f_i\left(x_i, x_{-i}, \xi_{ij}\right),
\end{equation}
where $p_i=[p_{ij}]_{j=1}^m\in \mathbb R^m$ is the probability vector belonging to the ambiguity set $P_i\subseteq \mathbb R^{m}$, and $\xi_{ij}$ are sampled scenarios. This formulation introduces a minimax structure, transforming each player’s problem into a saddle-point problem, where they minimize their objective under the worst-case probability distribution in the ambiguity set.
This {\it distributionally robust Nash equilibrium} (DRNE) framework ensures robust solutions, even under uncertainty. 

{\bf Contributions.} We study the problem of distributionally robust Nash equilibrium (DRNE), where multiple players interact in a noncooperative game under uncertainty. Although the stochastic Nash equilibrium problem \eqref{SNE} has been extensively analyzed, its distributionally robust counterpart remains mainly  underexplored, particularly in terms of efficient algorithms with provable convergence guarantees. Unlike standard minimax problems, which involve a single player solving a minimax objective, the DRNE problem involves $n$ players, each solving their own minimax problem. This distinction renders existing minimax algorithms unsuitable \cite{nedic2009subgradient,hamedani2021primal,zhao2022accelerated}, as they are designed for single-player settings. To address this gap, we make the following contributions:
\begin{enumerate}
\item We reformulate DRNE problems as \textit{Variational Inequalities (VI)}, providing a unified framework that enables the use of scalable and efficient solution methods.While \cite{pantazis2023data} discussed modeling DRNE as a variational inequality, in our work we consider nonsmooth and large-scale problems and we establish rigorous convergence guarantees.
\item We propose a gradient-based algorithm for solving DRNE problems, with two key theoretical results: (i) an oracle complexity of $\zz{{\mathcal{O}}\left(\frac{1}{\epsilon^{2}}\log^2(\frac{1}{\epsilon})\right)}$, and (ii) showing that the iterates generated by the proposed algorithm converge to a solution of the \eqref{DRO_Nash} problem almost surely. In contrast to \cite{juditsky2011solving}, which proposed a stochastic mirror-prox method for variational inequalities and established convergence in expectation, our analysis addresses set-valued monotone operators, defines the gap via a supremum over all subgradients to ensure well-posedness in nonsmooth settings, and further provides almost sure convergence guarantees.
\end{enumerate}
Next, we review related works on Variational inequalities and Distributionally Robust Nash Equilibrium (DRNE) problems, highlighting gaps in the existing literature. 
In Section \ref{reform VI}, we show that how the DRNE problem can be reformulated as a VI. In Section \ref{sec:pre}, we discuss some definitions and state our assumptions. In Section \ref{sec:alg}, we analyze a gradient descent-ascent method and establish its oracle complexity and convergence results. Finally, in Section \ref{sec:numeric}, we evaluate the performance of the proposed scheme through numerical experiments.
}

\subsection{Related Works}
\aj{Nash Equilibrium (NE) is a cornerstone of game theory, modeling optimal decision-making among self-interested agents \cite{nash1951non,nash1950equilibrium}. Nash established the existence of a mixed strategy NE for finite strategic games \cite{nash1950equilibrium}. Extensions to continuous strategy sets, where players choose from convex and compact sets, depend on specific assumptions about the continuity and convexity of the payoff functions \cite{bacsar1998dynamic}. In such settings, NE problems can be equivalently reformulated as Variational Inequality (VI) problems \cite{facchinei2003finite}, enabling the application of efficient computational methods. Over the past few decades, several methods have been proposed for solving VIs with guaranteed convergence rates \cite{korpelevich1976extragradient,doi:10.1137/S1052623403425629,malitsky2015projected,jalilzadeh2019proximal,alizadeh2024randomized}. Notable approaches include the Extragradient (EG) and Mirror-Prox (MP) methods \cite{juditsky2011solving,nemirovski2004prox}, the forward-backward-forward (FBF) algorithm \cite{tseng2000modified}, dual extrapolation \cite{nesterov2007dual}, and the reflected gradient/forward-reflected-backward (FoRB) scheme \cite{malitsky2020forward}. These algorithms achieve a complexity of $\mathcal O(1/\epsilon)$ for VIs with monotone and Lipschitz continuous operators in a deterministic setting. In recent years, stochastic VIs have gained prominence, particularly in applications involving uncertainty. Techniques such as Sample Average Approximation (SAA) and Stochastic Approximation (SA) have been employed to tackle these problems. SAA estimates the expectation of the stochastic mapping by averaging over a large set of samples, while SA updates the solution iteratively using mini-batch samples. For monotone VIs, the stochastic variants of these methods exhibit a complexity of $\mathcal O(1/\epsilon^2)$ \cite{nemirovski2009robust,juditsky2011solving}. Additional foundational results on the existence and uniqueness of Nash equilibria in stochastic settings can be found in \cite{ravat2011characterization} and \cite{jadamba2015variational}, who employed VI approaches to study stochastic Nash equilibria.

\textbf{Distributionally Robust Nash Equilibrium (DRNE).}
Distributionally Robust Optimization (DRO) has emerged as a practical approach for decision-making under uncertainty, balancing robustness and out-of-sample performance \cite{pantazis2023data}. In DRO, decision-makers optimize against the worst-case distribution within an ambiguity set, addressing uncertainties without the conservatism of robust optimization or the data-hungry nature of Sample Average Approximation (SAA). The extension of DRO to Nash games, known as Distributionally Robust Nash Equilibrium (DRNE), has gained attention in recent years. \cite{qu2012distributionally} introduced the concept of DRNE in finite games, providing foundational results. \cite{liu2018distributionally} studied its existence in continuous games, highlighting the challenges of developing numerical methods for DRNE. Other works, such as \cite{sun2016convergence}, have explored the stability properties of DRNE as players acquire more data. In \cite{pantazis2023data}, the authors study data-driven Wasserstein distributionally robust Nash equilibrium (DRNE) problems with heterogeneous uncertainty, where each agent builds a private ambiguity set from sample data. They \zz{ established variational inequality reformulations for this problem, }provide finite-sample guarantees and prove that the resulting equilibria converge almost surely to the true Nash equilibrium as data grows, highlighting robustness and consistency in learning equilibria from data. 
In contrast, our work reformulates the \eqref{DRO_Nash} with nonsmooth objective function as a variational inequality (VI) and develops a gradient descent–ascent algorithm with provable convergence guarantees.}

\section{Variational Inequality Reformulation of DRNE}\label{reform VI}
\aj{The problem of finding a Distributionally Robust Nash Equilibrium (DRNE) can be reformulated as a Variational Inequality (VI) problem, leveraging the minimax structure inherent in the DRNE framework. Unlike standard minimax problems, which involve a single player solving a minimax objective, the DRNE problem involves $n$ players, each solving their own minimax problem. This makes existing minimax algorithms unsuitable \cite{nedic2009subgradient,hamedani2021primal,zhao2022accelerated}, as they are designed for single-player setups. First, recall the DRNE problem defined earlier:
\begin{equation}\tag{DRNE}
\min _{x_i \in K_i} \max _{p_i \in P_i} \sum_{j=1}^m p_{ij} f_i\left(x_i, x_{-i}, \xi_{ij}\right),
\end{equation}
where $x_i$ and $x_{-i}$ represent the decisions of player $i$ its opponents, $P_i$ is the ambiguity set for the probability distribution, and $f_i$ is a convex but possibly nonsmooth function. To reformulate this as a VI, we use the optimality conditions for each player under the minimax structure. Define $F(x,p) \triangleq \begin{bmatrix}F_1(x,p) \\ F_2(x,p)\end{bmatrix}$, 
where $F_1$ and $F_2$ are the partial subdifferentials of the player objectives 
with respect to the primal and dual variables, respectively. 
Specifically, for each $i=1,\dots,n$ we set
$F_{1}(x,p)\triangleq \prod_{i=1}^n\left(\sum_{j=1}^m p_{ij}\,\partial_{x_i} f_i(x_i,x_{-i};\xi_{ij})\right),$ 
$F_{2}(x,p)\triangleq \left[\big[-f_i(x_i,x_{-i},\xi_{ij})\big]_{j=1}^m\right]_{i=1}^n$, $g_1(x,p)\in F_1(x,p)$, and $g_2(x,p)=F_2(x,p)$. 
Therefore, using first order optimality condition, solving problem \eqref{DRO_Nash} is equivalent to finding $z^* \in Z$ such that
$
\left\langle g(z^*), z-z^*\right\rangle \geqslant 0 \quad \forall z \in Z,
$
\za{w}here $g(z)\triangleq \begin{bmatrix}g_1(x,p) \\ g_2(x,p)\end{bmatrix}$, $z=\begin{bmatrix}x \\ p\end{bmatrix}$, $Z=K \times 
P$, $K=\prod_{i=1}^n K_i$ and $P=\prod_{i=1}^n P_i$. Note that the nonsmooth nature of $f_i$ implies that the operator $F$ in the VI formulation is also nonsmooth. This lack of smoothness requires the development of a novel algorithm, as most existing methods for solving VIs rely on smoothness assumptions. Despite this challenge, the VI reformulation simplifies the theoretical analysis and serves as a foundation for designing scalable algorithms to solve the DRNE problem efficiently. Next, we define a gap function to measure the quality of the solution obtained from the algorithm. 

\textbf{Gap Function.} 
Karamardian \cite{karamardian1976existence} showed that under continuity and pseudomonotonicity of a monotone operator, 
solving the variational inequality is equivalent to solving the Minty variational inequality (MVI) \cite{minty1962monotone}. 
In our nonsmooth setting, the operator $F$ is set--valued due to the use of subdifferentials. 
Accordingly, the Minty VI requires finding $z^* \in Z$ such that
$\langle g,\, z-z^*\rangle \;\geq\; 0 
\quad \forall z\in Z,\; \forall g \in F(z).$
The advantage of the Minty formulation is that it naturally admits a gap function representation, 
which can be used to quantify algorithmic progress \cite{alizadeh2024randomized}. 
In particular, since our operator is monotone, convergence guarantees can be established by analyzing the following gap function for any $z \in Z$:
$\mathrm{Gap}(z) \triangleq \sup_{y\in Z}\;\sup_{g\in F(y)} \langle g,\, z-y\rangle.$
}
\section{Preliminaries}\label{sec:pre}
In this section, we initially outline crucial notations,  followed by articulating the primary assumptions essential for the convergence analysis.

\textbf{Notations.} Throughout this paper, $\left\|.\right\|^2$ represents the Euclidean vector norm, $\mathcal{P}_{K}(u)$ is the projection of $u$ onto set ${K}$, i.e. $\mathcal{P}_K [x] = \mbox{argmin}_{y \in K} \| y-x \|$, and $\mathbb{E}[.]$ denotes the expectation of a random variable. Moreover, we denote by 
$\partial_{x_i} f_i(x)$ the subdifferential of $f_i$ with respect to $x_i$. 

\aj{Next, we state our assumptions on the problem structure, including the existence of an optimal solution and the properties of objective function in  \eqref{DRO_Nash} problem.}
\begin{assumption}\label{assump:F} \aj{(i) $K_i$ and $P_i$ are convex and compact sets for any $i=1,\hdots,n$.}
(ii) The set of optimal solutions is nonempty.
(iii) For each player $i=1,\ldots,n$, the function 
$f_i:K\to \mathbb{R}$ is convex in its own decision variable $x_i$ 
(possibly nonsmooth). 
The collection $\{f_i\}_{i=1}^n$ is assumed to satisfy the monotonicity condition
$\sum_{i=1}^n 
\big\langle g_i(x)-h_i(y),\,x_i-y_i\big\rangle \geq 0,$ for all $x,y\in K,
g_i(x)\in \partial_{x_i} f_i(x), h_i(y)\in \partial_{x_i} f_i(y).$
\end{assumption}

Note that (iii) in Assumption \ref{assump:F} guarantees that operator 
$F(z)$
is monotone on $K$.
\aj{Recall the definition of the partial subdifferential of the objective functions in problem \eqref{DRO_Nash}, we can write  \za{$g_1(x, p) =\Big[\sum_{j=1}^m p_{i j} \tilde{\nabla} f_i\left(x_i, x_{-i}; \xi_{i j}\right)\Big]_{i=1}^n,\quad
g_2(x, p) =\Big[\big[-f_i(x_{i,} x_{-i}, \xi_{i j})\big]_{j=1}^m\Big]_{i=1}^n,$}%
where $\tilde \nabla f_i\in \partial f_i$. Since $m\gg 1$, computing $g_1\left(x, p\right)$ and $g_2(x, p)$ exactly is computationally expensive. Instead, we estimate them using mini-batches. Specifically, we draw mini-batches $B_1, B_2$ from $\{1, \hdots , m\}$ uniformly at random without replacement, where $|B_1| = b_1 < m$ and $|B_2| = b_2 < m$.  
The mini-batch gradient estimates are defined as $g_{1, B_1}(x, p) = \left[g_{1, B_1}^{(i)}(x_i, x_{-i}, p_i)\right]_{i=1}^n\in F_{1,B_1}(x,p)$ and $g_{2, B_2}(x, p) = \left[g_{2, B_2}^{(i)}(x_i, x_{-i}, p_i)\right]_{i=1}^n\in F_{2,B_2}(x,p)$, where
$g_{1, B_1}^{(i)}(x_i, x_{-i}, p_i) \triangleq \frac{m}{b_1}\sum_{j \in B_1} p_{ij} \tilde{\nabla} f_i\left(x_i, x_{-i}; \xi_{ij}\right),\quad  g_{2, B_2}^{(i)}(x_i, x_{-i}, p_i) \triangleq \left[g_{2, B_2}^{(i), j}(x_i, x_{-i}, p_i)\right]_{j=1}^m,$
$$ g_{2, B_2}^{(i), j}(x_i, x_{-i}, p_i) \triangleq 
    \begin{cases} 
        0 & \text{if } j \notin B_2, \\
        -\frac{m}{b_2}f_i\left(x_i, x_{-i}; \xi_{ij}\right) & \text{if } j \in B_2,
    \end{cases}.$$}
Note that $g_{1,B_1}(x,p)$ and $g_{2,B_2}(x,p)$ are unbiased estimators of $g_1(x,p)$ are $g_2(x,p)$, i.e.,\break  $\mathbb{E}[g_{\ell,B_{\ell}}(x,p)-g_{\ell}(x,p)] =0,$  for all $x \in K$ and $p \in P$ and $\ell\in\{1,2\}$. Next, we make the following assumptions regarding the boundedness of their variance.
\aj{\begin{assumption}\label{assum:samples} There exists $\nu_\ell>0$ such that $\mathbb{E}[\|g_{\ell,B_{\ell}}(x,p)-g_{\ell}(x,p)\|^2] \leq \nu_\ell^2/b_{\ell}$ for all $g_{\ell,B_\ell}\in F_{\ell,B_\ell}(x,p)$, $g_\ell\in F_\ell(x,p)$, $x \in K$, $p \in P$ and $\ell\in\{1,2\}$.
\end{assumption}}
\aj{The next assumption ensures the boundedness of sample partial subgradients in expectation. } 
\begin{assumption}\label{ass:bounded_operator}
   \aj{There exist $M_x,M_p>0$ such that for any $g_{\ell,B_\ell}\in F_{\ell,B_\ell}(x,p)$}:  $\mathbb{E}\left[\left\|g_{1, B_1}\left(x, p\right)\right\|^2\right]\leqslant M_x^2$ and $\mathbb{E}\left[\left\|g_{2, B_2}\left(x, p\right)\right\|^2\right]\leqslant M_p^2,$ for all $  x \in {K}, p \in P.$
\end{assumption}\vspace{-3mm}
\aj{\begin{remark}
    Note that, from definition of $g_{1,B_1}$ and the fact that $P$ is a bounded set, one can show that $\mathbb{E}\left[\left\|g_{1, B_1}\left(x, p\right)\right\|^2\right]$ is bounded if subgradient of $f_i$s for all $i$ is bounded. Moreover, boundedness of $\mathbb{E}\left[\left\|g_{2, B_2}\left(x, p\right)\right\|^2\right]$ follows from the fact that $f_i$s are continious functions and $K$ is a bounded set.
\end{remark}}
\vspace{-5mm}
\section{Proposed Method}\label{sec:alg}
\aj{Now we are ready to propose the \textbf{Gradient Descent-Ascent for DRNE (GDA-DRNE)} algorithm, a stochastic optimization method designed to solve Distributionally Robust Nash Equilibrium (DRNE) problems. The details can be seen in Algorithm \ref{alg1}. The algorithm operates iteratively, alternating between updating the decision variables $x_i$ and the probability distributions $p_i$ for each player $i$. At each iteration $t$, mini-batches $B_1$ and $B_2$ are sampled randomly from the set of scenarios $\{1, \dots, m\}$, where 
$ |B_1| = b_1 < m$ and  $|B_2| = b_2 < m$,
to approximate the subgradient terms $g_{1, B_1}$ and $g_{2, B_2}$, respectively. 
The step sizes, $\lambda_t$ and  $\gamma_t$, are chosen from predefined positive sequences.}
\begin{algorithm}
	\caption{Gradient Descent-Ascent for DRNE (GDA-DRNE)} \label{alg1}
        \textbf{Input:}$x^0 \in X, p^0 \in P,\left\{\gamma_t, \lambda_t\right\}_{t \geq 0} \subseteq \mathbb{R}_{+};$

        \For{$t=0,\ldots T-1$}{
        Choose $B_1, B_2 \in \{1,\dots, m\}$ randomly such that $\left|B_1\right|=b_1<m$ and $\left|B_2\right|=b_2<m$\;
    		\For {$i=1,\dots n$}{
            Select $g_{1, B_1}\in F_{1,B_1}(x^t,p^t)$ and $g_{2, B_2}\in F_{2,B_2}(x^t,p^t)$,\vspace{-2mm}
 \begin{align*}
    & x^{t+1}_i\gets\mathcal{P}_{K_i} \left(x^{t}_i-\lambda_tg_{1, B_1}^{(i)}\left(x^t_i, x^t_{-i},p^t_i\right)\right);\\
    & p^{t+1}_i\gets\mathcal{P}_\ze{P_i}\left(p^t_i-\gamma_t g_{2, B_2}^{(i)}\left(x^{t}_i,x^{t}_{-i},p^t_i\right)\right);
    \end{align*}\vspace{-8mm}}
    	}	
\end{algorithm}\vspace{-5mm}
\subsection{Convergence Analysis}
In this section, we analyze the convergence properties of Algorithm \ref{alg1}. 
Proofs are in Appendix \ref{app:A}.
\begin{theorem}[Convergence rate]\label{Th: conv_rate}
    Suppose Assumptions \ref{assump:F},   \ref{assum:samples}, and \ref{ass:bounded_operator} hold. Let $\{x^t,p^t\}_{t}$ be the sequence generated by Algorithm \ref{alg1}, where $x=[x_i]_{i=1}^n$ and $p=[p_i]_{i=1}^n$. Choose the stepsizes as 
    $\lambda_t = \gamma_t=\frac{1}{\sqrt{1+t}\ log(t+2)}$ and define $\bar x^T\triangleq \tfrac{1}{ \sum_t\lambda_t}\sum_t\lambda_t x^t$ and $\bar p^T\triangleq \tfrac{1}{\sum_t\gamma_t}\sum_t\gamma_t p^t$, then \zz{for all $T \geq 1$} the following holds:  
\begin{align}\label{final rate}
    &\mathbb{E}\Big[\sup_{(x,p)\in K\times P}\zz{\sup_{\substack{g_1\in F_1(x,p) \\ g_2\in F_2(x,p)}}} \{\langle  g_1(x,p),\bar x^T-x\rangle+\langle \zz{g_2(x,p)},\bar p^T-p\rangle\}\Big]\leq \mathcal{O}\!\left(\tfrac{\log T}{\sqrt{T}}\right).
\end{align}
Consequently, the oracle complexity required to achieve an $\epsilon$-gap is 
    $\mathcal{O}\!\left(\tfrac{1}{\epsilon^{2}}\log^{2}\!\tfrac{1}{\epsilon}\right)$.
\end{theorem}
\begin{remark}
   Although the mapping with respect to $p$ is linear in our distributed formulation, implying that $F_2(x,p)$ is single-valued and $g_2 = F_2(x,p)$, we retain the supremum over $g_2$ in \eqref{final rate} to preserve generality. This notation ensures that the convergence bound remains valid even when $F_2$ is set-valued, i.e., in nonsmooth extensions of the problem.
\end{remark}
Next, we show that the
iterates generated by Algorithm \ref{alg1} converges to a solution of \eqref{DRO_Nash} problem almost surely. 
\begin{theorem}[Almost Sure Convergence]\label{Th:con_in_exp}
Consider iterates $\{x^t,p^t\}_{t}$  generated by Algorithm \ref{alg1} , where $x=[x_i]_{i=1}^n$ and $p=[p_i]_{i=1}^n$.
Suppose Assumptions \ref{assump:F},   \ref{assum:samples}, and \ref{ass:bounded_operator} hold and choose $\lambda_t = \gamma_t=\frac{1}{\sqrt{1+t}\ log(t+2)}$, then $\{x^t,p^t\}_{t}$ converges to a solution of  \eqref{DRO_Nash} problem almost surely.
\end{theorem}

\section{Numerical Experiments}\label{sec:numeric}

We consider a class of risk-averse Nash games (RaNE), where each player minimizes a convex risk measure of their random cost function. Several risk measures, including Conditional Value-at-Risk (CVaR), admit equivalent expectation-form representations \cite{rockafellar2002conditional}. For instance, CVaR at confidence level $\alpha$ can be written as
$\mathrm{CVaR}_\alpha(Z) = \inf_{u \in \mathbb{R}} \left\{ u + \tfrac{1}{1-\alpha} \mathbb{E}\big[(Z - u)_+\big] \right\}$ 
revealing an underlying stochastic expectation structure via the auxiliary variable $u$. 
In this experiment, we consider a $n$-player Nash game, where each player solves $\min_{x_i\in K_i,u_i\in \mathbb R} \mathbb E[\phi(h_i(x_i,x_{-i},\xi_i),u_i)]$, where $\phi(z,u) := u + \tfrac{1}{1-\alpha}(z-u)_+$ is convex but possibly nonsmooth.\vspace{-3mm} 
\begin{figure}[!htb]
	\centering
{{\includegraphics[scale=0.14]{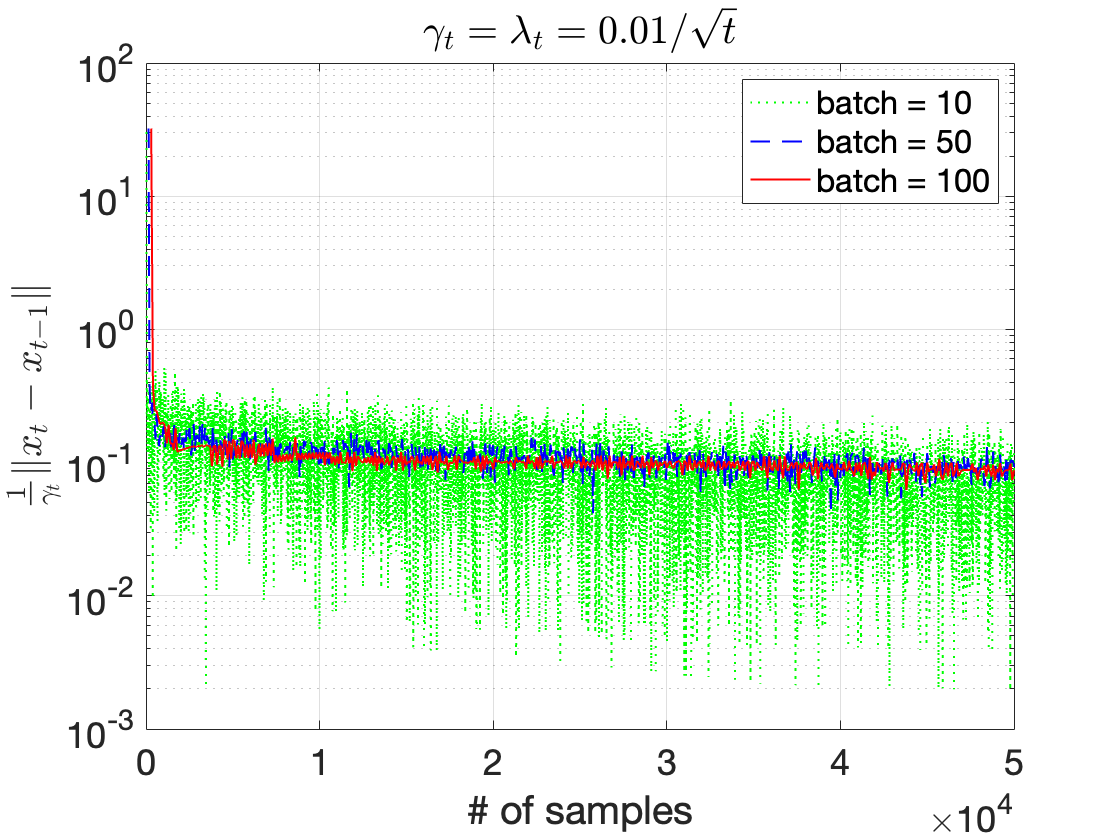} }}
   \hspace{0.3cm}
{{\includegraphics[scale=0.14]{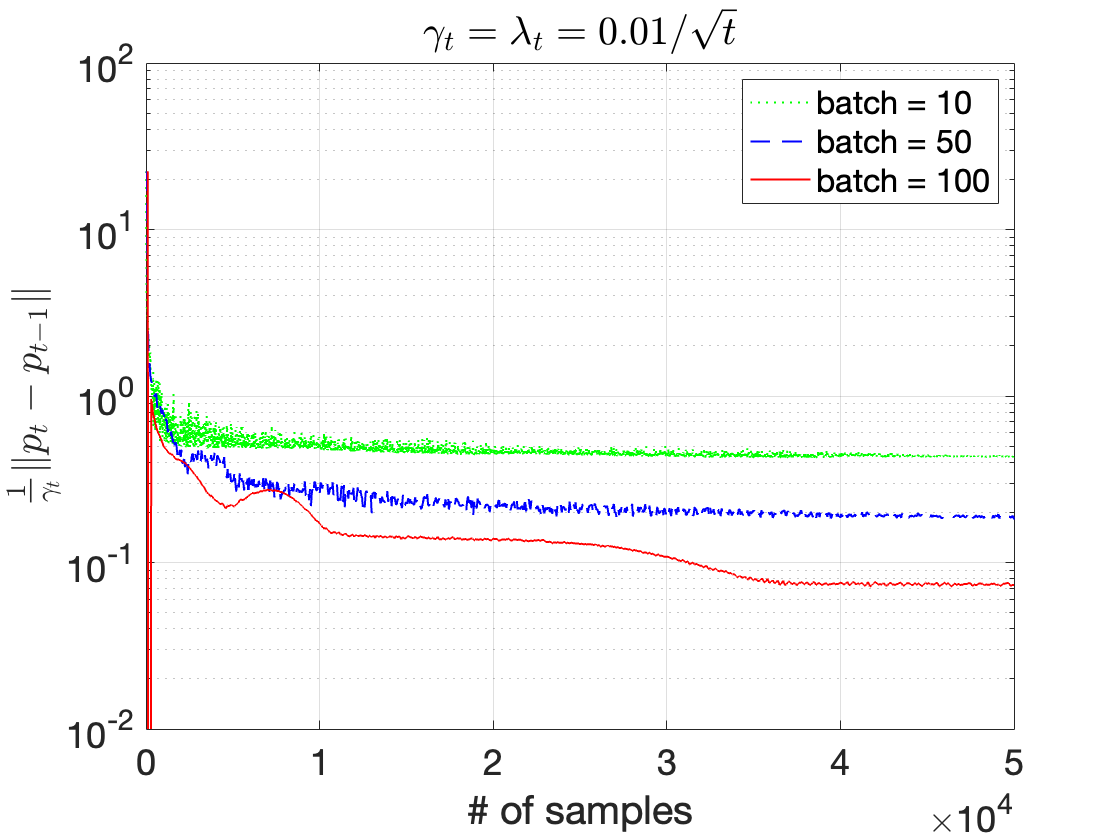} }}\vspace*{-3mm}
	\caption{ {Performance of GDA-DRNE averaged over 5 simulations for different batch sizes. }}
	\label{Fig:figurez4}%
\end{figure}
Following \eqref{DRO_Nash} setting, we consider the distributionally robust reformulation, where each player hedges against uncertainty in the probability distribution of $\xi$ by optimizing over a worst-case distribution within a simplex-defined ambiguity set $\Delta_i$, leading to $\min _{x_i \in K_i, u_i\in \mathbb R} \max _{p_i \in \Delta_i} \sum_{j=1}^m p_{ij} \phi(h_i\left(x_i,x_{-i}, \xi_{ij}),u_i\right)$. We implement GDA-DRNE on the risk-averse game with $n=5$ players, $n_i=10$, $m=100$, $\alpha=0.95$, $h_i(x,\xi_{ij})=\tfrac{1}{2}\xi^1_{ij}\|x\|^2+\xi^2_{ij}c^\top x$, where $x=(x_i,x_{-i})$, $\xi_{ij} = (\xi_{ij}^{(1)}, \xi_{ij}^{(2)})$ follow uniform distributions, $c$ follows standard normal distribution and $K_i=[-10,10]^{n_i}$. Figure \ref{Fig:figurez4} illustrates the performance of GDA-DRNE, demonstrating that the algorithm successfully converges to a solution.

\bibliography{biblio}
\appendix
\section{}\label{app:A}

The following lemma is essential for analyzing the convergence rate.

\begin{lemma}\label{lemma ineq} \za{(see Lemma 4 in \cite{alizadeh2024randomized})}
Given an arbitrary sequence $\{\bar{\sigma}_k\}_{k \geq 0} \subseteq \mathbb{R}^n$ and $\{\bar{\alpha}_k\}_{k \geq 0} \subseteq \mathbb{R}^{++}$, let $\{v_k\}_{k \geq 0}$ be a sequence such that $v_0 \in \mathbb{R}^n$ and $v_{k+1} = v_k + \frac{\bar{\sigma}_k}{\bar{\alpha}_k}$. Then, for all $k \geq 0$ and $x \in \mathbb{R}^n$,
\[
\langle \bar{\sigma}_k, x - v_k \rangle \leq \frac{\bar{\alpha}_k}{2} \|x - v_k\|^2 - \frac{\bar{\alpha}_k}{2} \|x - v_{k+1}\|^2 + \frac{1}{2\bar{\alpha}_k} \|\bar{\sigma}_k\|^2.
\]
\end{lemma}
Before proving the convergence result, we state the following definition:
\aj{\begin{definition}[Stochastic errors]\label{def:stoch_err} For our analysis, we define the following stochastic terms for $t\geq 0$:
 $w_1^t\triangleq g_{1, B_1}(x^t,p^t)-g_{1}(x^t,p^t)$ and  $w_2^t\triangleq g_{2, B_2}(x^t,p^t)-g_{2}(x^t,p^t)$.
\end{definition}}
\textbf{Proof of Theorem \ref{Th: conv_rate}.}
Defining $x=[x_i]_{i=1}^n$ and using the update rule of $x_i$ in Algorithm \ref{alg1},  \za{yields $x^{t+1}=\mathcal{P}_{K} \left(x^{t}-\lambda_t\zz{g_{1, B_1}\left(x^t, p^t\right)}\right).$ 
By} the definition of projection we have that: $$x^{t+1}=\arg \min _{x \in {K} }\left\|x-\left(x^t-\lambda_t \zz{g_{1, B_1}\left(x^t, p^t\right)}\right)\right\|^2.$$
Then, from the first-order optimality condition, one can obtain the following: $$\left\langle x^{t+1}-\left(x^t-\lambda_t  \zz{g_{1, B_1}\left(x^t, p^t\right)}\right), x -x^{t+1}\right\rangle \geqslant 0 \quad \forall x \in {K}.$$
Now rearranging the terms leads to:
\begin{align}\label{opt cond}
\nonumber \lambda_t\left\langle \zz{g_{1, B_1}\left(x^t, p^t\right)}, x^{t+1}-x\right\rangle  &\leqslant\left\langle x^{t+1}-x^t, x-x^{t+1}\right\rangle \\
& =\frac{1}{2}\left(\left\|x^t-x\right\|^2-\left\|x^{t+1}-x\right\|^2-\left\|x^{t+1}-x^t\right\|^2\right).
\end{align}
By adding and subtracting $x^t$ to the left hand side and using the fact that $ab\leq \tfrac{1}{2}a^2+\tfrac{1}{2}b^2$, we get:
\begin{align*} \lambda_t\langle \zz{g_{1, B_1}\left(x^t, p^t\right)}, x^t-x\rangle &\leqslant \frac{1}{2}\left(\left\|x^t-x\right\|^2-\left\|x^{t+1}-x\right\|^2-\left\|x^{t+1}-x^t\right\|^2\right)\\&\quad+\frac{\lambda_t^2}{2}\left\|\zz{g_{1, B_1}\left(x^t, p^t\right)}\right\|^2+\frac{1}{2}\left\|x^{t+1}-x^t\right\|^2\\
&=\frac{1}{2}\left(\left\|x^t-x\right\|^2-\left\|x^{t+1}-x\right\|^2\right)+\frac{\lambda_t^2}{2}\left\|\zz{g_{1, B_1}\left(x^t, p^t\right)}\right\|^2.
\end{align*}
Using the fact that \zz{$g_{1, B_1}\left(x^t, p^t\right)=g_{1}\left(x^t, p^t\right)+w_1^t$}, we obtain that:
\begin{align*}\lambda_t\langle  \zz{g_{1}\left(x^t, p^t\right)}, x^t-x\rangle  \leqslant\tfrac{1}{2}(\left\|x^t-x\right\|^2-\left\|x^{t+1}-x\right\|^2)+\tfrac{\lambda_t^2}{2}\left\|\zz{g_{1, B_1}\left(x^t, p^t\right)}\right\|^2+\lambda_t\langle w_1^t,x-x^t\rangle.
\end{align*}
Similarly, \ze{using the update rule of $p_i^{t+1}$} in Algorithm \ref{alg1} and defining $p=[p_i]_{i=1}^n$, one can obtain that:
\begin{align*} \gamma_t\langle \zz{g}_{2}\left(x^t, p^t\right), p^t-p\rangle  &\leqslant \tfrac{1}{2}\left(\left\|p^t-p\right\|^2-\left\|p^{t+1}-p\right\|^2\right)+\tfrac{\ze{\gamma}_t^2}{2}\left\|\zz{g}_{2, B_2}\left(x^t, p^t\right)\right\|^2+\gamma_t\langle w_2^t,p-p^t\rangle.
\end{align*}
Summing two inequalities and using the monotonicity of operator $F$, \zz{for any $g_1(x,p) \in F_1(x,p)$ and $g_2(x,p) \in F_2(x,p)$} we obtain that: 
\begin{align*} &\lambda_t\langle \zz{g_{1}\left(x, p\right), x^t-x\rangle}+\gamma_t\langle \zz{g_{2}\left(x, p\right)}, p^t-p\rangle \\&\quad \leqslant\frac{1}{2}\left(\left\|x^t-x\right\|^2-\left\|x^{t+1}-x\right\|^2\right)+\frac{\lambda_t^2}{2}\left\|\zz{g_{1, B_1}}\left(x^t, p^t\right)\right\|^2+\lambda_t\langle w_1^t,x-x^t\pm v^t\rangle\\
&\qquad+ \frac{1}{2}\left(\left\|p^t-p\right\|^2-\left\|p^{t+1}-p\right\|^2\right)+\frac{\ze{\gamma}_t^2}{2}\left\|\zz{g_{2, B_2}\left(x^t, p^t\right)}\right\|^2+\gamma_t\langle w_2^t,p-p^t\pm v^t\rangle.
\end{align*}
\za{Applying Lemma~\ref{lemma ineq} with $\bar{\alpha} = 1$ and $\bar{\sigma} = \lambda_t w_1^t$ yields a bound on $\langle \lambda_t w_1^t, x - v^t \rangle$, while using $\bar{\sigma} = \gamma_t w_2^t$ gives a bound on $\langle \gamma_t w_2^t, p - v^t \rangle$, \zz{for any $g_1(x,p) \in F_1(x,p)$ and $g_2(x,p) \in F_2(x,p)$} we obatin the following
}
\begin{align*} &\lambda_t\langle \zz{g}_{1}\left(x, p\right), x^t-x\rangle+\gamma_t\langle \zz{g}_{2}\left(x, p\right), p^t-p\rangle \\ &\leqslant \tfrac{1}{2}\left(\left\|x^t-x\right\|^2-\left\|x^{t+1}-x\right\|^2\right)+\tfrac{\lambda_t^2}{2}\left\|\zz{g}_{1, B_1}\left(x^t, p^t\right)\right\|^2+\lambda_t\langle w_1^t,v^t-x^t\rangle+\tfrac{1}{2}\|x-v^t\|^2\\&\quad-\tfrac{1}{2}\|x-v^{t+1}\|^2+\tfrac{\lambda_t^2}{2}\|w_1^t\|^2+ \tfrac{1}{2}\left(\left\|p^t-p\right\|^2-\left\|p^{t+1}-p\right\|^2\right)+\tfrac{\ze{\gamma}_t^2}{2}\left\|\zz{g}_{2, B_2}\left(x^t, p^t\right)\right\|^2\\&\quad+\gamma_t\langle w_2^t,v^t-p^t\rangle+\tfrac{1}{2}\|p-v^t\|^2-\tfrac{1}{2}\|p-v^{t+1}\|^2+\tfrac{\gamma_t^2}{2}\|w_2^t\|^2.
\end{align*}
Summing both sides \ze{over $t = 0,\dots,T-1$} and defining $\bar x^T\triangleq \tfrac{1}{ \sum_t\lambda_t}\sum_t\lambda_t x^t$ and $\bar p^T\triangleq \tfrac{1}{\sum_t\gamma_t}\sum_t\gamma_t p^t$, then taking the supremum over all selections $g_1(x,p) \in F_1(x,p)$ and $g_2(x,p) \in F_2(x,p)$ we get that:
\begin{align*}
    &\left(\sum_{t=0}^{T-1}\lambda_t\right)\zz{\sup_{g_1(x,p)}\langle  g_1(x,p),\bar x^T-x\rangle}+\left(\sum_{t=0}^{T-1}\gamma_t\right)\zz{\sup_{g_2\in F_2(x,p)}}\langle \zz{g}_2(x,p),\bar p^T-p\rangle\\
    &\quad \leq {1\over 2}\|x^0-x\|^2+{1\over 2}\|x-v^0\|^2+\sum_{t=0}^{T-1}{\lambda^2_t\over 2}\|\zz{g}_{1,B_1}(x^t,p^t)\|^2+\sum_{t=0}^{T-1}{\lambda^2_t\over 2}\|w_1^t\|^2+\sum_{t=0}^{T-1}{\lambda_t}\langle w_1^t,v^t-x^t\rangle\\
    &\qquad +{1\over 2}\|p^0-p\|^2+{1\over 2}\|p-v^0\|^2+\sum_{t=0}^{T-1}{\gamma^2_t\over 2}\|\zz{g}_{2,B_2}(x^t,p^t)\|^2+\sum_{t=0}^{T-1}{\gamma^2_t\over 2}\|w_2^t\|^2 +\sum_{t=0}^{T-1}{\gamma_t}\langle w_2^t,v^t-p^t\rangle.
\end{align*}

Taking the supremum over $(x,p)\in K\times P$, then taking the expectation from both sides of the aforementioned inequality, using \z{Assumptions} \ref{assum:samples}, and \ref{ass:bounded_operator},  considering $\gamma_t = \lambda_t $, and dividing both sides by $\sum_{t=0}^{T-1}\lambda_t$, implies that: 

\begin{align*}
    & \mathbb{E}[\sup_{(x,p)\in K\times P} \zz{\sup_{\substack{g_1\in F_1(x,p) \\ g_2\in F_2(x,p)}}} \{\langle  g_1(x,p),\bar x^T-x\rangle+\langle \zz{g_2(x,p)},\bar p^T-p\rangle\}]\\
    &\leq \tfrac{1}{\sum_{t=0}^{T-1}\lambda_t}\Big( \tfrac{1}{2}(B_X+B_P)+\sum_{t=0}^{T-1}\tfrac{\lambda^2_t}{2}\big(M_x^2+\tfrac{\nu_1^2}{b_1}+M_p^2+ \tfrac{\nu_2^2}{b_2}\big)\Big).
\end{align*}
Where $B_X\triangleq \sup_{x\in K}\{\|x-x^0\|^2+\|x-v^0\|^2\}$ and $B_P\triangleq \sup_{p\in P}\{\|p-p^0\|^2+\|p-v^0\|^2\}$. 
\af{By selecting $\lambda_t =\frac{1}{\sqrt{t+1}\log(t+2)}$, we observe that
\begin{align*}
\sum_{t=0}^{T}\lambda_t=\sum_{t=0}^{T}\frac{1}{\sqrt{t+1}\log(t+2)}
&\geq
\frac{1}{\log(T+2)}\sum_{t=0}^{T}\frac{1}{\sqrt{t+1}}\\&\geq
\frac{1}{\log(T+2)}\int_{0}^{T}\frac{dx}{\sqrt{x+1}}
=\frac{2(\sqrt{T+1}-1)}{\log(T+2)},
\end{align*}
and
\[
\sum_{t=0}^{T-1}\lambda_t^2 = \sum_{t=0}^{T} \frac{1}{(t+1)(\log(t+2))^2} \leq 4.
\]
Defining
\[
C_T \triangleq \frac{1}{2}(B_X + B_P) + 2(M_x^2+\tfrac{\nu_1^2}{b_1}+M_p^2+\tfrac{\nu_2^2}{b_2}),
\]
we can derive the following desired result for any $T> 1$:
\begin{align*}
    &\mathbb{E}\Big[\sup_{(x,p)\in K\times P}\zz{\sup_{\substack{g_1\in F_1(x,p) \\ g_2\in F_2(x,p)}}} \{\langle  g_1(x,p),\bar x^T-x\rangle+\langle \zz{g_2(x,p)},\bar p^T-p\rangle\}\Big]\\
    &\leq \frac{1}{\sum_{t=0}^{T-1}\lambda_t}\Big( \tfrac{1}{2}(B_X+B_P)+\tfrac{1}{2}\big(M_x^2+\tfrac{\nu_1^2}{b_1}+M_p^2+\tfrac{\nu_2^2}{b_2}\big)\sum_{t=0}^{T-1}\lambda_t^2 \Big)\\& \leq \frac{\log(T+2)C_T}{2\sqrt{T+1}-1}\leq \mathcal{O}\left(\tfrac{\log T}{\sqrt{T}}\right). \quad \square
\end{align*}}
To prove the almost sure convergence, we use the following technical lemma.

\begin{lemma} \label{lem1} \za{(Robbins-Siegmund)} Suppose $\omega_t,\eta_t,\upsilon_t
$ and $\psi_t$ are nonnegative random variables that satisfy:
$$
\mathbb{E}\left[\omega_{t+1} \mid \mathcal{F}_t\right] \leq\left(1+\eta_t\right) \omega_t+\upsilon_t-\psi_t, \quad \forall t\geq0.
$$
where $\mathbb{E}\left[\omega_{t+1} \mid \mathcal{F}_t\right]$ represents the conditional expectation for the given $\{\omega_0...,\omega_t\},\{\eta_0,...,\eta_t\},$\break$\{\upsilon_0,...,\upsilon_t\}
$, $\{\psi_0,...,\psi_t\}$, and $\sum_{t=0}^\infty \eta_t < \infty$, $\sum_{t=0}^\infty \upsilon_t < \infty$ , then we have:  $
\text { (i) } \lim _{t \rightarrow \infty} \omega_t$ exists and is finite, i.e. $\omega_t \rightarrow \omega$ almost surely, and $\text { (ii) } \sum_{t=0}^\infty \psi_t < \infty$.
    
\end{lemma}
\textbf{Proof of Theorem \ref{Th:con_in_exp}.}
\ze{By adding and subtracting $x^t$ to the left hand side of \eqref{opt cond} from the proof of Theorem \ref{Th: conv_rate}, and }using the fact that $ab\leq {1\over 2\alpha}a^2+{\alpha\over 2}b^2$ for some $\alpha\in(0,1)$ the following can be obtained:
\begin{align*} &\lambda_t\langle g_{1, B_1}\left(x^t, p^t\right), x^t-x\rangle \\&\leqslant \tfrac{1}{2}(\left\|x^t-x\right\|^2-\left\|x^{t+1}-x\right\|^2-\left\|x^{t+1}-x^t\right\|^2)+\frac{\lambda_t^2}{2\alpha}\left\|g_{1, B_1}\left(x^t, p^t\right)\right\|^2 +\tfrac{\alpha}{2}\left\|x^{t+1}-x^t\right\|^2.
\end{align*}
Using the fact that $g_{1, B_1}\left(x^t, p^t\right)=g_{1}\left(x^t, p^t\right)+w_1^t$, we obtain that:
\begin{align*} &\lambda_t\langle g_{1}\left(x^t, p^t\right), x^t-x\rangle \\&\leqslant\frac{1}{2}\left(\left\|x^t-x\right\|^2-\left\|x^{t+1}-x\right\|^2\right)+\frac{\alpha-1}{2}\left\|x^{t+1}-x^t\right\|^2+\frac{\lambda_t^2}{2\alpha}\left\|g_{1, B_1}\left(x^t, p^t\right)\right\|^2+\lambda_t\langle w_1^t,x-x^t\rangle.
\end{align*} 
Similarly, by using the update rule  $p^{t+1}=\mathcal P_\ze{P}\left(p^t-\gamma_t F_{2, B_2}\left(x^{t},p^t\right)\right)$ in Algorithm \ref{alg1} and following the same steps, we get:
\begin{align*} &\lambda_t\langle g_{2}\left(x^t, p^t\right), p^t-p\rangle \\&\leqslant\frac{1}{2}\left(\left\|p^t-p\right\|^2-\left\|p^{t+1}-p\right\|^2\right)+\frac{\alpha-1}{2}\left\|p^{t+1}-p^t\right\|^2+\tfrac{\gamma_t^2}{2\alpha}\left\|g_{2, B_2}\left(x^t, p^t\right)\right\|^2+\gamma_t\langle w_2^t,p-p^t\rangle.
\end{align*} 
Summing two inequalities and using the monotonicity of operator $F$, \zz{for any $g_1(x,p) \in F_1(x,p)$ and $g_2(x,p) \in F_2(x,p)$} we obtain that
\begin{align*} &\lambda_t\langle g_{1}\left(x, p\right), x^t-x\rangle+\gamma_t\langle g_{2}\left(x, p\right), p^t-p\rangle \\&\leqslant\frac{1}{2}\left(\left\|x^t-x\right\|^2-\left\|x^{t+1}-x\right\|^2\right)+\frac{\alpha-1}{2}\left\|x^{t+1}-x^t\right\|^2+\frac{\lambda_t^2}{2\alpha}\left\|F_{1, B_1}\left(x^t, p^t\right)\right\|^2+\lambda_t\langle w_1^t,x-x^t\rangle\\
&+\frac{1}{2}\left(\left\|p^t-p\right\|^2-\left\|p^{t+1}-p\right\|^2\right)+\frac{\alpha-1}{2}\left\|p^{t+1}-p^t\right\|^2+\tfrac{\gamma_t^2}{2\alpha}\left\|g_{2, B_2}\left(x^t, p^t\right)\right\|^2+\gamma_t\langle w_2^t,p-p^t\rangle.
\end{align*} 
For any solution $(x^*,p^*)$ of the VI reformulation of \eqref{DRO_Nash} and the fact that $\lambda_t,\gamma_t > 0$, there exist $g_1\in F_1(x^*,p^*)$ and $g_2\in F_2(x^*,p^*)$, such that $0\leq \lambda_t\langle g_{1}\left(x^*, p^*\right), x^t-x^*\rangle+\gamma_t\langle g_{2}\left(x^*, p^*\right), p^t-p^*\rangle $, hence
by taking the conditional expectation from both sides, using Assumptions \ref{assum:samples} and \ref{ass:bounded_operator} and unbiasedness of $w_1^t,w_2^t$ one can obtain: 
\begin{align*} &\nonumber\mathbb {E}\left[\|x^{t+1}-x^* \|^2| \mathcal{F}_t\right] + \mathbb E\left[\|p^{t+1}-p^* \|^2| \mathcal{F}_t\right]\\& \leqslant \|x^{t}-x^* \|^2 +\|p^{t}-p^* \|^2 +\frac{\lambda_t^2}{\alpha}M_x^2+\frac{\gamma_t^2}{\alpha}M_p^2-(1-\alpha)\big(\left\|x^{t+1}-x^t\right\|^2+\left\|p^{t+1}-p^t\right\|^2\big).
\end{align*} 
Choosing $\lambda_t= \gamma_t = \frac{1}{\sqrt{t+1}\log(t+2)}$ in the above inequality, the requirements of Lemma \ref{lem1} are satisfied, and we can conclude that:
$\{\left\|x^t-x^*\right\|^2+\left\|p^t-p^*\right\|^2\}_t$ is a convergent sequence almost surely, hence, $\{(x^t,p^t)\}_t$ is a convergent sequence almost surely, i.e., there exists $(x^\#,p^\#)\in K\times P$ such that $(x^t,p^t)\to (x^\#,p^\#)$ almost surely. 
Now from the convergence rate result in Theorem \ref{Th: conv_rate}, $$\lim_{T\to \infty} \mathbb{E}\Big[\sup_{(x,p)\in K\times P}\zz{\sup_{\substack{g_1\in F_1(x,p), g_2\in F_2(x,p)}}}  \{ \langle g_1(x,p),\bar x^T-x\rangle+\langle g_2(x,p),\bar p^T-p\rangle\}\Big]=0$$ and there exists $B\in(0,\infty)$ such that $$\sup_{(x,p)\in K\times P}\zz{\sup_{\substack{g_1\in F_1(x,p), g_2\in F_2(x,p)}}}  \{ \langle g_1(x,p),\bar x^T-x\rangle+\langle g_2(x,p),\bar p^T-p\rangle\}\leq B$$ due to compactness of $K\times P$, therefore, the dominated convergence theorem and continuity of the gap function imply that \za{$$\mathbb{E}\Big[\sup_{(x,p)\in K\times P}\zz{\sup_{\substack{g_1\in F_1(x,p), g_2\in F_2(x,p)}}}  \{\lim_{T\to \infty}( \langle g_1(x,p),\bar x^T-x\rangle+\langle g_2(x,p),\bar p^T-p\rangle)\}\Big]=0.$$}
Since $\sum_t \lambda_t=\sum_t \gamma_t =\infty$, then $(\bar x^T,\bar p^T)\to (x^\#,p^\#)$. Therefore, $$\mathbb{E}[\sup_{(x,p)\in K\times P}\zz{\sup_{\substack{g_1\in F_1(x,p),g_2\in F_2(x,p)}}}  \{ \langle g_1(x,p),x^\#-x\rangle+\langle g_2(x,p),p^\#-p\rangle\}]=0.$$ This implies that $ \langle g_1(x,p),x^\#-x\rangle+\langle g_2(x,p),p^\#-p\rangle=0$ almost surely for all $(x,p)\in K\times P$, $g_1\in F_1(x,p)$ and $g_2\in F_2(x,p)$. Which means that $(x^\#,p^\#)$ is a solution of \eqref{DRO_Nash} problem almost surely. $\square$
\end{document}